\input amstex
\documentstyle{amsppt}
\document
\magnification=1200
\NoBlackBoxes
\nologo
\vsize18cm


\bigskip

\centerline{\bf FUNCTIONAL EQUATIONS}
\smallskip

\centerline{\bf FOR QUANTUM THETA FUNCTIONS}

\medskip

\centerline{\bf Yuri I. Manin}

\medskip

\centerline{\it Max--Planck--Institut f\"ur Mathematik, Bonn}

\centerline{\it and Northwestern University, Evanston, USA}

\bigskip

{\it Abstract.} Quantum theta functions were introduced by the author 
in [Ma1]. They are certain elements in the function rings of quantum tori.
By definition, they satisfy a version of the classical functional 
equations involving shifts
by the multiplicative periods. This paper shows that for
a certain subclass of period lattices  (compatible with the
quantization form), quantum thetas satisfy
an analog of another classical functional
equation related to an action of the metaplectic group upon the (half of)
the period matrix. In the quantum case, this is replaced
by the action of the special orthogonal group on the quantization form,
which provides Morita equivalent tori. The argument
uses Rieffel's approach to the construction of (strong)
Morita equivalence bimodules and the associativity
of Rieffel's scalar products. 

\bigskip  

\centerline{\bf \S 0. Introduction and summary}

\medskip

{\bf 0.1. Theta functions and theta vectors.} This paper is a contribution
to the theory of quantum theta functions introduced in
[Ma1] and further studied in [Ma2], [Ma3]. It addresses two
interrelated questions:

\smallskip

(a) {\it What is the connection between quantum theta
functions and theta vectors?}  

This question was repeatedly raised by
A.~S.~Schwarz, see e.~g. [Sch].

\smallskip

(b) {\it Does there exist a quantum analog of 
the classical functional equation for thetas (related to the action
of the metaplectic group, see e.~g. [Mu], \S 8)?}

\smallskip

Briefly, the (partial) answers we give here look as follows.

\smallskip

(i) Schwarz's theta vectors are certain elements
of projective modules over $C^{\infty}$-- or $C^*$-- rings
of unitary quantum tori. When such a module
is induced from the basic Heisenberg representation 
by a lattice embedding into a vector Heisenberg group,
the respective theta vectors $f_T$ are parametrized by
the points $T$ of Siegel upper half space, and in different
models of the basic representation take the form
of a ``quadratic exponent'' $e^{\pi ix^tTx}$, a classical theta,
or Fock's vacuum state: see Theorem 2.2 in [Mu].

\smallskip

To the contrary, quantum thetas are certain elements 
of the $C^{\infty}$ function ring itself. (For this reason, partial multiplication of quantum thetas studied in [Ma3] does not seem
to be directly related to the tensor product of projective (bi)modules).

\smallskip

The basic relationship between the two classes of objects is this.
For a lattice embedded in a vector Heisenberg group,
Rieffel's scalar products of {\it theta vectors} (these products take values in
the $C^{\infty}$ ring of the relevant quantum torus)
are certain {\it quantum theta functions}. This extends a 
calculation of Section 3 in
[Ma4], which in turn generalized a result of [Bo]:
see Theorems 3.2.1 and 3.6.1 below. Theorem 3.5.1 characterizes
in abstract terms the subclass of quantum thetas that
can be obtained in this way.

\smallskip

(ii) The classical functional equation relates two thetas
considered as sections of line bundles
over two isomorphic complex tori (Fourier series).
Bundles and sections are lifted 
to the universal covers which are then identified
compatibly with period lattices. 

\smallskip

Similarly, the 
functional equation for scalar product quantum thetas stated here
relates two theta functions in two quantum tori algebras
related by a bimodule inducing
their  Morita equivalence. The equation then simply says
that the respective thetas coincide after being applied to appropriate
vacuum vectors, and becomes a particular case of Rieffel's associativity relations:
see the Theorem 3.3.1.

\smallskip

Here are some details. Consider a classical theta function of $z\in \bold{C}^N$
$$
\theta (z,\Omega ):=\sum_{n\in \bold{Z}^N} e^{\pi i\, n^t\Omega n +2\pi i\, n^tz}
\eqno(0.1)
$$
where $\Omega$ is a symmetric complex matrix with positive defined
imaginary part. This function satisfies two sets of functional equations.
Firstly, for all $m\in \bold{Z}^N$,
$$
\theta (z+m,\Omega )=\theta (z,\Omega ), 
\eqno(0.2)
$$
$$
\theta (z+\Omega m,\Omega )=
e^{-\pi i\,m^t\Omega m -2\pi i\,m^tz}\,\theta (z,\Omega ) .
\eqno(0.3)
$$
Secondly,
$$
\theta (\Omega^{-1}z,-\Omega^{-1} )= 
(\roman{det}\, (\Omega/i))^{1/2}\,e^{\pi i\,z^t\Omega^{-1}z}\,\theta (z,\Omega ) .
\eqno(0.4)
$$
In fact, (0.4) is the most important special case of a more general modular
functional equation related to the action of $Sp (2,\bold{Z})$ upon
the space of pairs $(z, \Omega )$ which we do not spell out here.

\smallskip

The geometric meaning of these equations can be described
as follows. Consider
$\theta (z,\Omega )$ as a global section of the trivial
line bundle over $\bold{C}^N$. Equations (0.2) allow us to
consider it as a global section of the trivial
line bundle over  $(\bold{C^*})^N$ as well. This section
is written as a Laurent series
in the basic characters $e(n):=e^{2\pi i\,n^tz}$ of $(C^*)^N$. Equations
(0.3) allow us to descend one step further, now
turning $\theta$ into a section of a nontrivial line
bundle $\Cal{L}$ on the complex torus $\bold{C}^N/D$
where $D$ is the sublattice generated by the unit vectors
and the columns of $\Omega$. This is achieved by
embedding $D$ into a vector Heisenberg group acting
upon $\bold{C}^N\times \bold{C}$ compatibly with the
projection, and then taking the quotient of this space
with respect to $D$: see e.~g. [Mu], p.~35. 
Now, changing the initial basis of  $\bold{C}^N$ (e.~g. 
replacing it by  the columns of $\Omega$)
produces an isomorphic triple $(\Cal{T},\Cal{L},\theta )$
consisting of a complex torus,
line bundle, and its section. This is the source of equation
(0.4) and its generalizations.

\smallskip

Now deform the multiplication rule of the characters $e(n)$
by choosing an antisymmetric real matrix $A$ and putting
$e_{\alpha}(m)\,e_{\alpha}(n) := \alpha (m,n)\,e_{\alpha} (m+n)$
where $\alpha (m,n)=e^{2\pi i\,m^tAn}$. The deformed characters
generate various function rings representing the quantum torus
$T(\bold{Z}^N,\alpha )$ which should be considered
as a deformation of $(\bold{C^*})^N$ or of its maximal compact subtorus.
The series (0.1) in which $e^{2\pi i\,n^tz}$ is replaced by
$e_{\alpha}(n)$ furnishes an example of quantum thetas,
 studied in [Ma1] -- [Ma4].
There, especially in [Ma3], a theory of the functional equations
of the type (0.2)--(0.3) is developed,
applicable to the quantum tori over $p$--adic fields as well.

\smallskip

In this paper, we propose an analog of the functional equation (0.4)
corresponding this time to the change of the quantization
matrix $A\mapsto A^{-1}.$ The noncommutative geometric
context replacing the classical isomorphism of triples 
$(\Cal{T},\Cal{L},\theta )$ invoked above, involves now
the (strong) Morita equivalence of the relevant quantum tori,
compatible complex structures on these quantum tori, and
theta vectors in the respective projective bimodule:
see [PoS] and the references therein.
The whole emerging picture is surprisingly parallel
to the classical one.

\medskip 

{\bf 0.2. Plan of the paper.} In \S 1, we recall the basic definitions
related to various Heisenberg groups we use in this paper,
sketch their representation theory, and reproduce the description
of quantum thetas given in [Ma3]. In \S 2,
we recall how Heisenberg representations produce projective
modules over quantum tori via lattice embeddings,
and sum up the main properties of Rieffel's
scalar products in this context. In \S 3, we elaborate
and prove the statements in (i), (ii) above for
vector Heisenberg groups and their extensions by finite groups.
The last subsection sketches some suggestions for further research.

\bigskip

\centerline{\bf \S 1. Heisenberg groups and their representations}

\medskip

{\bf 1.1. Central extensions.} Let $\Cal{K}$ (resp. $\Cal{Z}$)
be an abelian group written additively (resp. multiplicatively).
Consider a function $\psi :\,\Cal{K}\times \Cal{K}\to \Cal{Z}$.
Then the following conditions (a) and (b) are equivalent:

\smallskip

(a) {\it $\psi (0,0)=1$ and $\psi$ is a cocycle, that is, for each $x,y,z$
we have}
$$
\psi (x,y)\psi (x+y,z) =\psi (x,y+z)\psi (y,z) .
\eqno(1.1)
$$

(b) {\it The following composition law on $\Cal{G}:=\Cal{Z}\times\Cal{K}$
turns $\Cal{G}=\Cal{G} (\Cal{K},\psi)$ into a group with identity $(1,0)$:}
$$
(\lambda ,x)(\mu ,y):= (\lambda\mu\psi (x,y), x+y).
\eqno(1.2)
$$

Moreover, if (a), (b) are satisfied, then the maps
$\Cal{Z}\to\Cal{G}:\,\lambda\mapsto (\lambda ,0),$
$\Cal{G}\to \Cal{K}:\, (\lambda ,x)\mapsto x$, 
describe $\Cal{G}$ as a central extension of
$\Cal{K}$ by $\Cal{Z}$:
$$
1\to \Cal{Z}\to \Cal{G}(\Cal{K},\psi )\to \Cal{K}\to 1.
\eqno(1.3)
$$
Notice that any bicharacter $\psi$ automatically satisfies (a).
For arbitrary $\psi$, putting $x=0$ in (1.1), we see that
$\psi (0,y)=1$ so that 
$$
(\lambda ,x) =(\lambda ,0)(1,x).
$$

\smallskip

{\bf 1.1.1. Bicharacter $\varepsilon$.} Consider any central extension
(1.3), choose a set theoretic section $\Cal{K}\to \Cal{G}:\,x\mapsto
\tilde{x}$ and define the map $\varepsilon :\,
\Cal{K}\times \Cal{K}\to \Cal{Z}$ by
$$
\varepsilon (x,y):= \tilde{x}\tilde{y}\tilde{x}^{-1}\tilde{y}^{-1}.
\eqno(1.4)
$$
Then $\varepsilon$ is a bicharacter which does not depend on
the choice of a section and which is antisymmetric: $\varepsilon (y,x)=
\varepsilon (x,y)^{-1}$, $\varepsilon(x,x) =1.$ In particular,
if $K\subset \Cal{K}$ is a subgroup liftable to $\Cal{G}$,
then $K$ is $\varepsilon$--isotropic.
\smallskip

For the group
$\Cal{G} (\Cal{K}, \psi )$, choosing $\tilde {x}= (1,x)$,
we find
$$
\varepsilon (x,y) =\frac{\psi (x,y)}{\psi (y,x)},
\eqno(1.5)
$$
and if $\psi$ itself is an antisymmetric bicharacter, then
$\varepsilon (x,y) =\psi (x,y)^2.$ 

\medskip

{\bf 1.1.2. Example: Heisenberg groups of quantum tori and quantum
theta functions.} 
Let $H$ be a free abelian
group of finite rank written additively, $k$ a ground field, and $\alpha :\,H\times H\to 
k^*$ a skewsymmetric pairing.
The quantum torus $T(H,\alpha )$
with the character group $D$ and quantization
parameter $\alpha$ is represented by an algebra
generated by a family of formal exponents $e(h)=e_{H,\alpha}(h),\,h\in H$,
satisfying the relations
$$
e(g) e(h)= \alpha (g,h) e(g+h).
\eqno(1.6)
$$
In particular, $T(H,1)$ is an algebraic torus, spectrum
of the group algebra $k[e(h)\,|\, h\in H]$ of $H$. The group of its
points $x\in T(H,1)(k) = \roman{Hom}\, (H,k^*)$ acts upon
functions on $T(H,\alpha )$ mapping $e_{H,\alpha}(h)$
to $x^*(e_{H,\alpha}(h)):=h(x) e_{H,\alpha}(h)$ where $h(x)$
denotes the value of the character $e(h)$ at $x$. 
 
\smallskip

The Heisenberg group  of $T(H,\alpha )$ introduced
in [Ma3] and denoted there $\Cal{G} (H,\alpha )$  consists of   
all maps of the form
$$
\Phi \mapsto c\, e_{H,\alpha}(g)\, x^*(\Phi )\, e_{H,\alpha}(h)^{-1},\
c\in k^*;\, x\in T(H,1)(k);\, g,h\in H.
\eqno(1.7)
$$
Any such map has a unique representative of the same form
in which $h=0$ (``left representative''). Writing this
representative as $[c; x, g]$ we get the composition law
$$
[c^{\prime};x^{\prime},g^{\prime}][c;x,g]=
[c^{\prime}c\, g(x^{\prime})\, \alpha(g^{\prime},g);x^{\prime}x,g^{\prime}+g].
\eqno(1.8)
$$
In other words, this group is the central extension of
$\roman{Hom}\,(H,k^*)\times H$ by $k^*$ corresponding
to the bicharacter
$$
\psi ((x^{\prime},g^{\prime}),(x,g))=g(x^{\prime}) \alpha (g^{\prime},g)
\eqno(1.9)
$$
and having the associated bicharacter
$$
\varepsilon ((x^{\prime},g^{\prime}),(x,g))=g(x^{\prime})
g^{\prime}(x)^{-1} \alpha^2(g^{\prime},g).
\eqno(1.10)
$$
In particular, if a subgroup $B\subset \roman{Hom}\,(H,k^*)\times H$
is liftable to $\Cal{G} (H,\alpha )$, the form (1.9)
restricted to $B$ must be symmetric: this is the main part
of Lemma 2.2 in [Ma3]. 

\smallskip

A lift $\Cal{L}$ of $B$ to a subgroup of $\Cal{G} (H,\alpha )$ is called
{\it a multiplier}. The  restriction to $B$ of the form (1.9),
$\langle \,,\rangle :\, B\times B\to k^*$, is called {\it the structure
form} of this multiplier. (Formal) linear combinations of the
exponents $e_{H,\alpha}$ invariant with respect to the
action of $\Cal{L} (B)$ constitute a linear space $\Gamma (\Cal{L})$
and are called (formal) {\it quantum theta functions}. 

\medskip

{\bf 1.2. Representations.} Given $\Cal{K},\,\Cal{Z},\, \psi$
and a ground field $k$ as above, choose in addition a character
$\chi :\, \Cal{Z} \to k^*$. Consider a linear space of functions
$f:\,\Cal{K}\to k$ invariant with respect to the affine shifts 
and define operators $U_{(\lambda ,x)}$
on this space by
$$
(U_{(\lambda ,x)}f) (x):= \chi (\lambda \psi(x,y)) f(x+y).
\eqno(1.11)
$$ 
A straightforward check shows that this is a representation of
$\Cal{G} (\Cal{K}, \psi )$. However, it is generally reducible.
Namely, suppose that there is an $\varepsilon$--isotropic subgroup
$K_0 \subset \Cal{K}$ liftable to $\Cal{G} (\Cal{K}, \psi ).$ Let
$\sigma :\, K_0\to \Cal{G} (\Cal{K}, \psi ),$ $\sigma (y) = (\gamma (y), y)$
be such a lift. Denote by $F(\Cal{K}//K_0)$ the subspace
of functions satisfying the following condition:
$$
\forall\, x\in \Cal{K},\, y\in K_0,\
(U_{(\gamma (y) ,y)}f) (x):= \chi (\varepsilon (x,y)) f(x),
\eqno(1.12)
$$  
or, equivalently,
$$
\forall\, x\in \Cal{K},\, y\in K_0,\
f(x+y) =\chi (\gamma (y)^{-1} \psi (y,x)^{-1}) f(x).
\eqno (1.13)
$$
This subspace is invariant with respect to (1.11).

\smallskip

Formula (1.13) shows that if we know the value
of $f$ at a point $x_0$ of $\Cal{K}$, it extends
uniquely to the whole coset $x_0+K_0$, hence the notation
$F(\Cal{K}//K_0)$ suggesting ``twisted'' functions
on the coset space $\Cal{K}/K_0$.

\smallskip

Clearly, a minimal subspace if this kind is obtained  
if we choose for $K_0$ a maximal isotropic subgroup.

\medskip

{\bf 1.3. Locally compact abelian topological groups.}
The formalism briefly explained above is only
an algebraic skeleton. In the category of $LCAb$ of
locally compact abelian topological groups
and continuous homomorphisms,
with properly adjusted definitions, one can get
a much more satisfying picture.

\smallskip

First af all, choose 
$\Cal{Z} := \bold{C}_1^*= \{ z\in \bold{C}^*\,|\, |z|=1\}$.
This is {\it a dualizing object}: for each $\Cal{K}$
in $LCAb$ there exists the internal $\roman{Hom}\, (\Cal{K},\Cal{Z})$
object, called the character group $\widehat{\Cal{K}}$,
and the map $\Cal{K} \mapsto \widehat{\Cal{K}}$ extends
to the equivalence of categories $LCAb \to LCAb^{op}$
(Pontryagin's duality).

\smallskip

Let now $\psi$ be a continuous cocycle $\Cal{K}\times\Cal{K}\to\Cal{Z}$
so that $\varepsilon$ is a continuous bicharacter.
Call the extension $\Cal{G} (\Cal{K}, \psi )$
a  {\it Heisenberg group}, if the map 
$x\mapsto \varepsilon (x,*)$ identifies $\Cal{K}$ with
$\widehat{\Cal{K}}$. 

\smallskip

Choose $k=\bold{C}$, and $\chi$ continuous. The formula (1.11) 
makes sense 
e.g. for continuous functions $f$.
Especially interesting, however, is the representation
on $L_2(\Cal{K})$ which makes sense  because the operators (1.11)
are unitary with respect to the squared norm $\int_{\Cal{K}} |f|^2
d\mu_{Haar}$. Of course, square integrable functions
cannot be evaluated at points, so that $f(x+y)$ in (1.11) must
be understood as the result of shifting $f$ by $y\in \Cal{K}$;
similar precautions should be taken in the formula (1.13)
defining now the space $L_2(\Cal{K}//K_0)$ where $K_0$
is a closed isotropic subgroup (it is then automatically
liftable to a closed subgroup),
and in many intermediate calculations. See Mumford's treatment on
pp. 5--11 of [Mu] specially tailored for
readers with algebraic geometric sensibilities.

\smallskip

The central fact of the representation theory of a
Heisenberg group $\Cal{G} (\Cal{K},\psi )$,
$\Cal{K}\in LCAb$, $\chi = id$, is this:

\medskip

{\it (i) If $K_0$ is a maximal isotropic subgroup,
$L_2(\Cal{K}//K_0)$ is irreducible.
 
\smallskip

(ii) Any unitary representation of 
 $\Cal{G} (\Cal{K},\psi )$ whose restriction
to the center is the multiplication by the identical character
is isomorphic to
the completed tensor product of $L_2(\Cal{K}//K_0)$
and a trivial representation. In particular,
representations upon $L_2(\Cal{K}//K_0)$
corresponding to different choices of $K_0$
are isomorphic.}
 
\medskip

For example, if $\Cal{K} = K_0\times \widehat{K}_0$
is a direct product of two maximal isotropic real
spaces, it has also maximal isotropic subgroups which are
sublattices in $\Cal{K}$, and the respective models
of the Heisenberg representation are connected by
a non--trivial isomorphism.

\medskip

{\bf 1.4. Variants and complements.} The category $LCAb$  
offers a clear--cut case of the representation theory of Heisenberg groups
whose further axiomatization seems elusive.
Nevertheless, the following general features of this case
persist in one form or another when one replaces
abelian groups by group schemes, objects of an abelian category etc.:

\medskip

{\it (A) $\Cal{Z}$ must be a dualizing object; Heisenberg
groups are singled out among other central extensions
by the condition that $\varepsilon$ identifies
$\Cal{K}$ with the dual object.

\smallskip

(B) A Heisenberg group has an essentially unique
representation upon twisted functions on
$\Cal{K}/K_0$ where $K_0$ is a maximal 
$\varepsilon$--isotropic subgroup.}

\medskip

Even in $LCAb$ and for the case of a real space $\Cal{K}$,
a meaningful and important variation of the principle
(B) occurs, when we allow to replace $K_0$
by a maximal isotropic subspace in the complexification
of $\Cal{K}$ (see [Mu]). We will recall and use
the respective construction of the Fock space model
in \S 3 below.   

\smallskip

It might be interesting to work out a similar formalism
in a DG and derived setting. For example,
in the category of abelian algebraic groups
$\bold{G}_m$ has many properties of a honest
dualizing group, however the dual object for
an abelian  variety $A$ is $\roman{Ext}^1 (A, \bold{G}_m)$
rather than $\roman{Hom}\, (A, \bold{G}_m).$

\bigskip

\centerline{\bf \S 2. Quantum tori and projective modules}

\medskip

{\bf 2.1. Embedded lattices and tori.} In this section $\Cal{K}$ denotes 
an object of $LCAb$, $\psi$ is a
bicharacter of $\Cal{K}$ such that 
$\varepsilon$ (cf. (1.5)) identifies $\Cal{K}$ with
$\widehat{\Cal{K}}$. Let $\Cal{G} (\Cal{K},\psi )$ be
the respective Heisenberg group, central extension of
$\Cal{K}$ by $\Cal{Z} = \bold{C}_1^*$ as above.

\smallskip

We will call {\it an embedded lattice} 
a closed subgroup $D\subset \Cal{K}$ such that
$D$ is a finitely generated free abelian group,
whereas $\Cal{K}/D$ is a topological torus, i.e.
a finite product of $S^1$. In this section we
consider only those groups $\Cal{K}$ which
admit embedded lattices. 

\smallskip

Consider a family of constants $c_h\in \bold{C}_1^*$,
$h\in D$. Put 
$$
E(h):= (c_h,h)\in \Cal{G} (\Cal{K},\psi )
\eqno(2.1)
$$ 
From (1.2) we get
$$
E(g)E(h)=\frac{c_gc_h}{c_{g+h}}\,\psi (g,h)\,E(g+h).
$$
Assume that
$$
\alpha (g,h) := \frac{c_gc_h}{c_{g+h}}\,\psi (g,h)
\eqno(2.2)
$$
is a skewsymmetric pairing. Then the map $e_{D,\alpha}(h)
\mapsto E(h)$ is compatible with the relations (1.6),
and in particular any representation $U$ of
$\Cal{G} (\Cal{K},\psi )$ induces a representation of
an appropriate function algebra of the quantum torus
$T(H,\alpha )$. One easily sees that any $\alpha$ on $D$ can be induced
from an appropriate lattice embedding of $D$;
one can even take $\psi$ to be a skewsymmetric bicharacter so that
$\alpha$ will coincide with the restriction of $\psi$.

\smallskip
We will consider two function
algebras of the quantum torus $T(D,\alpha ).$
The algebra $C^{\infty}(D,\alpha )$ of smooth functions consists
of infinite series $\sum_{h\in D}a_he_{D,\alpha}(h)$
where the formal exponents satisfy (1.6), and
coefficients $\{ a_h\in \bold{C}\,|\,h\in D\}$
belong to the Schwarz's space $S(D)$. This algebra is
endowed with involution $(\sum_{h\in D}a_he_{D,\alpha}(h))^*=
\sum_{h\in D}\overline{a}_he_{D,\alpha}(h)^{-1}$. 

\smallskip

The $C^*$--algebra $C^{*}(D,\alpha )$ can be defined
as the universal algebra generated by the unitaries
$e_{D,\alpha}(h)$ satisfying (1.6). More concretely,
consider the action of $C^{\infty}(D,\alpha )$
upon $L_2(D)$ which is given by extending the multiplication
in $C^{\infty}(D,\alpha )$. Complete $C^{\infty}(D,\alpha )$
with respect to the operator norm. The result
will be $C^{*}(D,\alpha )$.

\smallskip 

Alternatively, any
element of $C^{*}(D,\alpha )$ can also be written
as a formal series  $\sum_{h\in D}a_he_{D,\alpha}(h)$
but there is no transparent way to specify
which sequences $\{ a_h\in \bold{C}\,|\,h\in D\}$
can occur as their ``noncommutative Fourier coefficients''.

\smallskip

A warning: the reader should not mix $\Cal{G} (\Cal{K},\psi )$
with the Heisenberg group of $T(D,\alpha )$ invoked in 1.1.2:
these two groups have totally different structures.

\medskip

{\bf 2.2. Inducing the Heisenberg representation.} Let
$(\lambda ,x)\mapsto U_{(\lambda ,x)}$ be an irreducible
unitary representation of $\Cal{G} (\Cal{K},\psi )$
in a Hilbert space  $\Cal{H}$. Since $\Cal{K}$ admits
an embedded lattice $D$, it is an abelian Lie group
(not necessarily connected) whose Lie algebra can be
identified with the tangent space to $\Cal{K}/D$ at zero.
The Heisenberg group $\Cal{G} (\Cal{K},\psi )$ is a Lie group as well.
Let $L$ be its Lie algebra.
A vector $f\in \Cal{H}$ is called {\it smooth} if
for any $X_1, \dots ,X_n \in L$ the following expression
makes sense
$$
\delta U_{X_1}\circ \dots \delta U_{X_n} (f)
$$
where $\delta U_{X} (f)$ is defined as the limit when $t\to 0$
$$
\delta U_X(f) := \roman{lim}\,\frac{U_{exp(tX)}f - f}{t}.
\eqno(2.3)
$$
It is known that the space $\Cal{H}_{\infty}$ of smooth vectors 
is dense and the operators $\delta U_X$ are skew adjoint but unbounded.

\medskip

{\bf 2.2.1. Theorem.} {\it The map $e_{D,\alpha}(h) \mapsto
E(h)$ induces on $\Cal{H}_{\infty}$ the structure of a 
finitely generated projective left $C^{\infty}(D,\alpha )$--module.}

\medskip

This module has an additional structure: scalar product with values in
$C^{\infty}(D,\alpha )$. Namely, assume that the scalar product
$\langle\,,\rangle$ in $\Cal{H}$ is antilinear in the second argument,
and put for $\Phi , \Psi\in \Cal{H}_{\infty}$
$$
{}_D\langle \Phi ,\Psi\rangle :=
\sum_{h\in D} \langle \Phi , e_{D, \alpha}(h) \Psi\rangle\,e_{D,\alpha}(h).
\eqno(2.4)
$$
Then this formal sum lies in $C^{\infty}(D,\alpha )$ and 
has the following properties:

\medskip

(i) {\it Symmetry:} ${}_D\langle \Phi ,\Psi\rangle^*=
{}_D\langle \Psi ,\Phi\rangle .$

\smallskip

(ii) {\it (Bi)linearity:} ${}_D\langle a\Phi ,\Psi\rangle =
a\,{}_D\langle \Phi ,\Psi\rangle$ for any $a\in C^{\infty}(D,\alpha )$.

\smallskip

(iii) {\it Positivity:}  ${}_D\langle \Phi ,\Phi\rangle$ 
belongs to the cone of positive elements of $C^{\infty}(D,\alpha )$.
Moreover, if ${}_D\langle \Phi ,\Phi\rangle =0$ then $\Phi =0.$

\smallskip

(iv) {\it Density:} The image of $\langle\,,\rangle$ is dense
in $C^{\infty}(D,\alpha )$.  

\medskip

{\bf 2.2.2. Theorem.} {\it The completion of $\Cal{H}_{\infty}$
with respect to the norm $\|\Phi\|^2:= \|{}_D\langle \Phi ,\Phi\rangle\|$
(where the rhs means the norm in $C^{\infty}(D,\alpha )$)
is a finitely generated projective left $C^{*}(D,\alpha )$--module $P$.
The scalar product (2.4) has a natural extension to
${}_D\langle\,,\rangle :\, P\times P \to C^{*}(D,\alpha )$.
The properties (i) -- (iv) hold for this extension as well.}

\medskip

{\bf 2.3. Dual embedded lattices.} Let $D\subset \Cal{K}$ be an embedded
lattice as in 2.1. Denote by $D^!$ the maximal closed subgroup
of $\Cal{K}$ orthogonal to $D$ with respect to $\varepsilon$.
From the Pontryagin duality it follows that $D^!$ (resp. $D$) can be
canonically identified with the character group of 
$\Cal{K}/D$ (resp. $\Cal{K}/D^!$) so that $D^!$ is an embedded lattice as well.

\smallskip

Assume moreover that $\psi$ is an antisymmetric pairing, so that one
can choose $E(h) = (1,h)\in \Cal{G} (\Cal{K}, \psi )$ for $h\in D$
and for $h\in D^!$ and 
define on $\Cal{H}_{\infty}$ the structure of $C^{\infty}(D^!,\alpha^!)$--module
as well where $\alpha^!$ is the pairing induced
on $D^!$ by $\psi$. Operators $e_{D,\alpha}(h),\,h\in D,$ commute 
with operators $e_{D^!,\alpha^!}(g),\,g\in D^!.$ 

\smallskip

We can consider $\Cal{H}_{\infty}$ as a right 
$C^{\infty} (D^!,\alpha^!)^{op}$--module. Moreover, we can and will 
identify the latter algebra with $C^{\infty} (D^!,\overline{\alpha}^!)$
by $e_{D^!,\alpha^!}(h) \mapsto e_{D^!,\overline{\alpha}^!}(h)^{-1}$
and extending this map by linearity.

\smallskip 

{\bf 2.3.1. Theorem.} {\it (i) We have 
$\|{}_D\langle \Phi ,\Phi\rangle\|^{1/2}
=\|{}_{D^!}\langle \Phi ,\Phi\rangle\|^{1/2}$
The completion $\Cal{H}$ of $\Cal{H}_{\infty}$
with respect to this norm is a projective left
module over both tori $C^{*}(D,\alpha )$ and $C^{*}(D^!,\alpha^!)$,
and each of these algebras is a total commutator of the other one.

\smallskip

(ii) Let  $C^{*}(D^!,\overline{\alpha}^!)$ 
act upon $\Cal{H}$ on the right as explained above. Consider the analog
of the scalar product (2.4)
$$
\langle \Phi ,\Psi\rangle_{D^!} :=\frac{1}{\roman{vol}\,\Cal{K}/D}\,
\sum_{h\in D} \langle  e_{D^!, \alpha^!}(h) \Psi ,\Phi 
\rangle\,e_{D^!,\overline{\alpha}^!}(h) \, \in C^{*}(D^!,\overline{\alpha}^!)
\eqno(2.5)
$$
It satisfies relations similar to (i)--(iv), and moreover,
for any $\Phi ,\Psi , \Xi$ the following
associativity relation holds:
$$
{}_D\langle \Phi ,\Psi\rangle \Xi = \Phi\,\langle \Psi ,\Xi\rangle_{D^!} .
\eqno(2.6)
$$
}
\smallskip

For proofs and further generalizations, see Rieffel's paper [Ri5],
in particular, sections 2 and 3.

\bigskip
\centerline{\bf \S 3. K\"ahler structure, theta--vectors, and quantum thetas} 

\medskip

{\bf 3.1. Case of vector Heisenberg groups.} 
In the first half of this section, we consider the case 
$\Cal{K} =$ a real vector space, $\psi$ an antisymmetric
bicharacter with values in $\bold{C}^*_1.$

\smallskip

In this case $\psi$ can be written in the form 
$$
\psi (x,y) = e^{\pi i A(x,y)}
\eqno(3.1)
$$
where $A:\, \Cal{K}\times \Cal{K} \to \bold{R}$
is a nondegenerate antisymmetric pairing. Choosing an appropriate
basis, we can identify $\Cal{K}$ with the space of pairs
of column vectors $x=(x_1,x_2),\, x_i\in \bold{R}^N$,
such that
$$
A(x,y) = x_1^ty_2-x_2^ty_1
$$ 
where $x_i^t$ denotes the transposed row vector.
We have then $\varepsilon (x,y) =e^{2\pi iA(x,y)}$.
In particular, the subspace $x_2=0$ is a maximal
$\varepsilon$--isotropic closed subgroup. Similarly, 
$\bold{Z}^{2N}$ is
a maximal $\varepsilon$--isotropic embedded
lattice.

\smallskip

We will recall the structure of two Heisenberg
representations of $\Cal{G} (\Cal{K},\psi )$ 
using normalizations adopted in [Mu].

\medskip

{\it Model I.} It consists of square integrable
complex functions $f$ on $\bold{R}^N$,
the first half of $\Cal{K}$, with the scalar product
$$
\langle f,g\rangle := \int f(x_1)\overline{g(x_1)}\,d\mu_x
\eqno(3.2)
$$ 
where $d\mu_x$ is the Haar measure in which $\bold{Z}^N$
has covolume 1.

\smallskip

The action of $\Cal{G} (\Cal{K},\psi )$, with central character
$\chi (\lambda ) =\lambda$, is given by the formula
$$
(U_{(\lambda ,y)}f) (x_1)= \lambda\,
e^{2\pi i\,x_1^ty_2 +\pi i\,y_1^ty_2}\,f(x_1+y_1).
\eqno(3.3)
$$

\smallskip
{\it Model II${}_T$.} The second model is actually
a family of models depending on the choice 
of a {\it compatible} K\"ahler structure
upon $\Cal{K}.$
A general K\"ahler structure on $\Cal{K}$
can be given by a pair consisting of a complex structure
and an Hermitean scalar product $H$. We will call this
K\"ahler structure compatible (with
the choice (3.1)) if $\roman{Im}\,H=A$. Such structures
are parametrized by the Siegel space consisting
of symmetric matrices $T\in M(N,\bold{C})$ with
positive defined $\roman{Im}\,T.$ 

\smallskip

In particular, the complex structure defined by $T$
identifies $\bold{R}^{2N}$ with $\bold{C}^N$
via 
$$
(x_1,x_2)=x \mapsto \underline{x} := Tx_1+x_2,
\eqno(3.4)
$$
and we have
$$
H(\underline{x},\underline{x})=\underline{x}^t\,(\roman{Im}\,T)^{-1}\,
\underline{x}^*
\eqno(3.5)
$$
where $*$ denotes the componentwise complex conjugation.

\smallskip

Consider the Hilbert space $\Cal{H}_T$ of holomorphic functions
on $\bold{C}^N=\Cal{K}$ consisting of the functions with finite
norm with respect to the scalar product
$$
\langle f, g\rangle_T := \int_{\bold{C}^g} f(\underline{x})\,
\overline{g(\underline{x})}\,
e^{-\pi H(\underline{x},\underline{x})}\,d\nu
\eqno(3.6)
$$
where $d\nu$ is the translation invariant measure making
$\bold{Z}^{2N}$ a lattice of covolume 1 in $\bold{R}^{2N}.$

\smallskip

For $(\lambda ,y)\in \Cal{G}(\Cal{K},\psi )$ and a holomorphic
function $f$ on $\Cal{K}$, put
$$
(U^{\prime}_{(\lambda ,y)}f)(\underline{x}):=
\lambda^{-1}e^{-\pi H(\underline{x},\underline{y})
-\frac{\pi}{2} H(\underline{y},\underline{y}) }
f(\underline{x}+\underline{y}).
\eqno (3.7)
$$
A straightforward check shows that these operators are 
unitary with respect to 
(3.6), and moreover, that they define a representation
of $\Cal{G}(\Cal{K},\psi )$ in $\Cal{H}_T$
corresponding to the character 
$\chi(\lambda )=\lambda^{-1}$ of $\bold{C}_1^*$,
in the sense of formula (1.11). This is (a version of) the classical
Fock representation.

\smallskip

It turns out that this representation is irreducible
and thus is an (antidual) model of the Heisenberg representation.

\smallskip

The proof of irreducibility spelled out in [Mu] 
involves constructing {\it vacuum vectors}
in $\Cal{H}_T$ which in this model turn out to be simply
constant functions. Translated via
canonical (antilinear) isomorphism into other models they look
differently, for example they become
(proportional to) a ``quadratic exponent''
$f_T:=e^{\pi ix^t_1 Tx_1}$ in Model I ( i.~e. $L_2(\bold{R}^{2N}//\bold{R}^{N})$) 
or to an essentially classical theta--function
$e^{\pi ix_1^t \underline{x}}\vartheta (\underline{x},T)$
in $L_2(\bold{R}^{2N}//\bold{Z}^{2N})$.
They are called  ``theta--vectors'' in [Sch].
For details, see the Theorem 2.2 in [Mu] and the discussion
around it. 

\medskip

{\bf 3.2. Scalar product ${}_D\langle *,*\rangle$.} We will now use
Model I
in order to induce projective modules
over toric algebras corresponding
to lattice embeddings $D\subset \Cal{K}$ as in 2.2 above.
Since our $\psi$ is already antisymmetric, we may 
and will put $c_h=1$ in (2.1), so that $\alpha$ is the restriction
of $\psi$ to $D$.
The main result of this subsection is the following
calculation.

\medskip

{\bf 3.2.1. Theorem.} {\it (i) We have
$$
{}_D\langle f_T, f_T\rangle =\frac{1}{\sqrt{2^N\,\roman{det\,Im}\,T}}
\sum_{h\in D} e^{-\frac{\pi}{2}\,H(\underline{h},\underline{h})}\,
e_{D,\alpha}(h).
\eqno(3.8)
$$
Moreover, 
$$
\Theta_D:=\sum_{h\in D} e^{-\frac{\pi}{2}\,H(\underline{h},\underline{h})}\,
e_{D,\alpha}(h)
$$ 
is a quantum theta function 
in the ring $C^{\infty} (D,\alpha )$ satisfying
the following functional equations:
$$
\forall\, g\in D,\ C_ge_{D,\alpha}(g)x_g^*(\Theta_D)=\Theta_D,
\eqno(3.9)
$$
where
$$
C_g=e^{-\frac{\pi}{2}\,H(\underline{g},\underline{g})},\ 
x_g^*(e_{D,\alpha}(h)) = e^{-\pi\,H(\underline{g},\underline{h})}\,
e_{D,\alpha}(h).
\eqno(3.10)
$$

(ii) We have 
$$
\Theta_D\,\bold{1}=
\sum_{h\in D} 
e^{-\pi H(\underline{h},\underline{h})-\pi H(\underline{x},\underline{h})} 
\eqno(3.11)
$$
where $\bold{1}$ is the vacuum vector in the model II${}_T$ represented by
the function identically equal to 1.}
\smallskip

{\bf Remarks.} (i) In the language  of [Ma3]
and 1.1.2 above, (3.8) means that $\Theta_D$
is invariant with respect to the multiplier $D\subset \Cal{G} (D,\alpha )$
where $D$ is embedded in the Heisenberg group
of our torus via  $g\mapsto [C_g; x_g, g]$ (left representatives).

\smallskip

(ii) The function  $\Theta_D \bold{1}$ is {\it complex conjugate} 
to the classical theta function
corresponding to a principal polarization of the complex torus
$\bold{C}^g/D$. Notice that this complex torus
is embedded into (the space of points of) the algebraic torus $T(D,1)(\bold{C}) = \roman{Hom}\,(D,\bold{C}^*)$
as its compact subtorus $\roman{Hom}\,(D,\bold{C}_1^*)$. 

\medskip

{\bf Proof of Theorem 3.2.1.} (i) We have, using (2.4), 
$$
{}_D\langle f_T, f_T\rangle =
\sum_{h\in D} \langle f_T, U_{(1,h)}f_T\rangle\,
e_{D,\alpha}(h) .
$$
From (3.3) we find
$$
(U_{(1 ,h)}f_T)(x_1) =
e^{\pi i\, (x_1^t+h_1^t)T(x_1+h_1)
+ 2\pi i\,x_1^th_2 +\pi i\,h_1^th_2} .
\eqno (3.12)
$$
so that in view of (3.2) 
$$
\langle f_T, U_{(1,h)}f_T\rangle =
\int_{\bold{R}^N}\, 
e^{\pi i\,x_1^tTx_1 -\pi i\, (x_1^t+h_1^t)\overline{T}(x_1+h_1)
- 2\pi i\,x_1^th_2 -\pi i\,h_1^th_2}\,d\mu.
\eqno(3.13)
$$
It remains to calculate the Gaussian 
integral in (3.13). 

\smallskip
 The exponent under the integral sign 
can be represented as $e^{-\pi\,(q(x_1)+l_h(x_1)+c_h)}$ where
$$
q(x_1)=2\,x_1^t\,\roman{Im}\,T\,x_1,\ 
l_h(x_1)=2 i \,x_1^t\,(\overline{T}h_1+h_2),\,\
c_h= ih_1^t(\overline{T} h_1 +h_2).
\eqno(3.14)
$$
Notice that $Th_1+h_2=\underline{h}$ in the notation (3.4),
so that $\overline{T}h_1+h_2=\underline{h}^*$ where $*$
denotes the componentwise complex conjugation.

\smallskip
We can solve for $\lambda_h\in \bold{C}^N$ the equation
$$
q(x_1+\lambda_h)-q(\lambda_h) =q(x_1)+l_h(x_1)
$$
We get
$$
\lambda_h=\frac{i}{2}\,(\roman{Im}\,T)^{-1} \underline{h}^*.
\eqno(3.15)
$$
Therefore
$$
\int e^{-\pi (q(x_1)+l_h(x_1)+c_h)}d\mu=
e^{-\pi (c_h-q(\lambda_h))}\int e^{-\pi q(x_1+\lambda_h)} d\mu=
\frac{e^{-\pi (c_h-q(\lambda_h))}}{\sqrt{\roman{det}\,q }}.
\eqno(3.16)
$$
Directly calculating $c_h-q(\lambda_h)$ we get 
$$
\frac{1}{2}\,\underline{h}^t (\roman{Im}\,T)^{-1}\,\underline{h}^* = 
\frac{1}{2}\,H(\underline{h},\underline{h}).
$$
Moreover, $\roman{det}\,q = 2^N \roman{det\, Im}\,T$
in view of (3.14). Hence we
finally recover (3.8). 

\smallskip

Now from (3.10)  we deduce
$$
C_ge_{D,\alpha}(g)\,x_g^*(\sum_{h\in D} e^{-\frac{\pi}{2}\,H(\underline{h},\underline{h})}\,
e_{D,\alpha}(h))=
$$
$$
e^{-\frac{\pi}{2}H(\underline{g},\underline{g})}\, e_{D,\alpha}(g)
\, \sum_{h\in D} 
e^{-\frac{\pi}{2}\,H(\underline{h},\underline{h})
-\pi H(\underline{g},\underline{h})}\,
e_{D,\alpha}(h)=
\sum_{h\in D} 
e^{-\frac{\pi}{2}\,H(\underline{g}+\underline{h},
\underline{g}+\underline{h})}\,
e_{D,\alpha}(g+h)
$$
because $e_{D,\alpha}(g)\,e_{D,\alpha}(h)=
e^{\pi i\,\roman{Im}\,H(\underline{g},\underline{h})}\,e_{D,\alpha}(g+h).$
This establishes (3.9).

\smallskip

(ii) Formula (3.11) now follows from (3.7) and (3.8).

\medskip

{\bf 3.3. A functional equation for quantum thetas: the basic case.}
Denote now by $D^!$ the dual embedded lattice as in 2.3.
From (3.1) one easily deduces that
$$
D^! =\{ x\in \Cal{K}\,|\,\forall \, y\in D,\ A(x,y)\in \bold{Z}\} .
\eqno (3.17)
$$
In particular, $D^!=\roman{Hom}\,(D,\bold{Z})$.

\smallskip

The following theorem is in fact a particular case of the
associativity formula (2.6) written for $\Phi = \Psi = \Xi = f_T$.
We replace one $f_T$ by the vacuum vector
$\bold{1}$ in the Model II${}_T$ and sketch an independent 
proof because this explicitly shows 
the structure of our functional equation.

\medskip

{\bf 3.3.1. Theorem.} {\it We have the following functional
equation for the pair of quantum theta functions 
$$
\Theta_D\bold{1} = \Theta_{D^!}\bold{1} .
\eqno(3.18)
$$
In other words,
$$
\sum_{h\in D}
e^{-\pi H(\underline{h},\underline{h})-\pi H(\underline{x},\underline{h})}=
\sum_{g\in D^!}
e^{-\pi H(\underline{g},\underline{g})-\pi H(\underline{x},\underline{g})}
\eqno(3.19)
$$
as functions of $x\in \Cal{K}.$
}

\smallskip

{\bf Proof.} (3.19) will follow from the Poisson summation formula
if we check the following. Put
$$
f_x(h):=e^{-\pi H(\underline{h},\underline{h})-\pi H(\underline{x},\underline{h})}
$$
considered {\it as a function of $h\in \Cal{K}$} ($x$ now being a parameter).
Define its Fourier transform by
$$
\widehat{f}_x(g) =\int_{\Cal{K}} f_x(h)\,e^{-2\pi iA(g,h)} d\nu_h.
$$
Then in fact  $\widehat{f}_x =f_x$.

\smallskip

The argument is similar to that in the proof of  Theorem 3.2.1. 
Denote by $Q(x)$ the real positive quadratic form $H(\underline{x},\underline{x})$
on $\Cal{K}$ considered as a real space. Put $R= \roman{Re}\,T,\,
S= \roman{Im}\,T$. After a somewhat tedious but 
straightforward calculation we find the matrix of $Q(x)$
written in real coordinates $(x_1,x_2)$:
$$
Q(x)=x_1^t(RS^{-1}R+S)\,x_1 + 2\,x_1^tRS^{-1}\,x_2 +
x_2^t S^{-1}x_2 =
$$
$$
(x_1^t,x_2^t)\, \left(
\matrix 
RS^{-1}R+S & RS^{-1}\\
S^{-1}R & S^{-1}\\
\endmatrix \right)\,
\left( \matrix 
x_1\\
x_2\\
\endmatrix \right)
$$
 
Now, find $\eta = (\eta_1 ,\eta_2 )\in \bold{C}\otimes_{\bold{R}}\Cal{K}$ from
the equation
$$
Q(h+\eta )- Q(\eta )= 
H(\underline{h},\underline{h}) + H(\underline{x},\underline{h}) 
+ 2i A(g,h).
\eqno(3.20)
$$
The result is
$$
\eta_1 = \frac{1}{2}\,x_1 -\frac{i}{2}\,
[(S^{-1}R (x_1+g_1) +S^{-1} (x_2+g_2)],
$$
$$
\eta_2 = \frac{1}{2}\,x_2 +\frac{i}{2}\,
[(S+RS^{-1}R) (x_1+g_1) +RS^{-1}(x_2+g_2)].
\eqno(3.21)
$$
Then we calculate $Q(\eta )$:
$$
Q(\eta )= -H(\underline{g},\underline{g}) - H(\underline{x},\underline{g}).
\eqno(3.22)
$$
Finally, from (3.20) and (3.22) we find
$$
\widehat{f}_x(g) =\int_{\Cal{K}} e^{-\pi q(h+\eta) +\pi q(\eta )} d\nu_h =
e^{\pi q(\eta )} = f_x(g).
$$

\medskip

{\bf 3.4. Comment: comparison of $\alpha$ and $\alpha^!$.} Choose a basis $(h_k\,|\,k=1,\dots ,N)$
of $D$ and consider the skew--symmetric matrix $A:= (A(h_k,h_l))$.
It determines the quantization parameter $\alpha$ of $T(D,\alpha )$.
To calculate $\overline{\alpha}^!$, choose a dual basis $(g_l)$ of $D^!$
determined by the condition $A(g_k,h_l)=\delta_{kl}$. A straightforward
calculation shows that 
$$
A^!:=A(g_k,g_l) = -A^{-1}
$$
so that finally $\overline{\alpha}^! (g_k,g_l)=e^{\pi i\,(A^{-1})_{kl}}.$

\smallskip

The map $A\mapsto A^{-1}$ is one of the standard generators of the group
$O\,(N,N; \bold{Z})/(\pm 1)$ acting on the strong Morita equivalence
classes of the quantum tori: cf. [RiSch] and [Li].

\medskip

{\bf 3.5. Comment: invariant characterization of 
quantum theta functions of the form $\Theta_D$.} In the language
of [Ma3], \S 2, $\Theta_D$ is the generator
of the space $\Gamma (\Cal{L})$ where $\Cal{L}$ is
the ample multiplier $\Cal{L} :\,D\to \Cal{G} (D, \alpha )$,
$\Cal{L} (g) = [C_g; x_g, g]$ (left representatives),
$C_g$ and $x_g$ being defined by (3.9).  

\smallskip

In fact, we have to check the conditions of [Ma3], Theorem 2.4.1.
The structure bilinear form of $\Cal{L}$ ( see [Ma3], (2.2)) is
$$
\langle h_1,h_2\rangle = h_2(x_{h_1})\,\alpha (h_1,h_2) =
=e^{-\pi\,H(\underline{h}_1,\underline{h}_2)}\,
e^{\pi i\,H(\underline{h}_1,\underline{h}_2)} =
e^{-\pi\,\roman{Re}\,H(\underline{h}_1,\underline{h}_2)} .
\eqno(3.23)
$$
Since $\roman{log}\, |\langle h,h\rangle | = -H(\underline{h},\underline{h})$
is a negative defined quadratic form, and the projection
$D\to \Cal{G} (D,\alpha )\to D$ is the (identical) isomorphism,
$\Cal{L}$ is ample and $\roman{dim}\,\Gamma (\Cal{L})=1.$

\smallskip

Moreover, as this calculation shows, $\Cal{L}$ has the following additional properties:

\smallskip

(a) The structure bilinear form of $\Cal{L}$ is real.

\smallskip

(b) There exists a K\"ahler structure upon $\bold{R}\otimes D$
consisting of a complex structure and an Hermitean form $H$
such that 
$$
\langle g,h\rangle = e^{-\pi \roman{Re}\,H(\underline{g},\underline{h})},\ 
\alpha (g,h) = e^{\pi i\,\roman{Im}\,H(\underline{g},\underline{h})}
\eqno(3.24)
$$ 
for all $g,h\in D$.

\smallskip

A converse statement is also true.

\medskip

{\bf 3.5.1. Theorem.} {\it Let $T(D,\alpha )$ be a quantum torus over $\bold{C}$,
with unitary quantization form $\alpha$. Let 
$\Cal{L}:\,B \to \Cal{G} (D,\alpha )$
be an ample multiplier such that the left representative projection
$B\to D$ (denoted also $h^{-}$ in [Ma3]) is an isomorphism
which we will use to identify $B$ with $D$.

\smallskip

Assume moreover that  one can define a K\"ahler structure on 
$\bold{R}\otimes D$ such that (3.24) holds ( $\langle\,,\rangle$
being the structure form of $\Cal{L}$).

\smallskip

Then there exists a real space $\Cal{K}$ endowed with
a bicharacter $\psi$, an compatible K\"ahler structure,
and a lattice embedding $D\subset \Cal{K}$
such that $\psi$ induces $\alpha$ on $D$,
and an appropriate generator of $\Gamma (\Cal{L})$ is of the form
$\Theta_D$ as above.
}
\medskip

{\bf Proof.} To see this, one should simply reverse the arguments above.
Take $\Cal{K} = \bold{R}\otimes D$ with the tautological embedding
of $D$, choose the K\"ahler structure such that (3.24) holds and
calculate the coefficients of an arbitrary generator
of $\Gamma (\Cal{L})$  as in [Ma3], (2.7). We will get
the right hand side of (3.8), up to a multiplicative
constant.
 
\medskip

{\bf 3.6. A generalization.} In this subsection, we will generalize
the Theorem 3.2.1 to the case of a lattice embedding
$D\subset \bold{R}^{2N}\times F\times \widehat{F}$
where $F$ is a finite group. 

\smallskip

Define the bicharacter $\psi_0$ of $F\times \widehat{F}$
by 
$$
\psi_0((a,l), (a^{\prime}, l^{\prime})):= l^{\prime}(a) \, .
\eqno(3.25)
$$
Via projection, we may and will consider it as a bicharacter
on $\bold{R}^{2N}\times F\times \widehat{F}$.
Similarly, $\psi$ from (3.1) induces
a bicharacter of $\bold{R}^{2N}\times F\times \widehat{F}$
which we denote by the same letter.
Consider the Heisenberg group $\Cal{G} (\bold{R}^{2N}\times F\times \widehat{F},
\psi \psi_0)$.  The map
$\Cal{G} (\bold{R}^{2N},
\psi )\times \Cal{G}( F\times \widehat{F},\psi_0) \to
\Cal{G} (\bold{R}^{2N}\times F\times \widehat{F},
\psi \psi_0)$,
$$
((\lambda , x), (\mu , a, l))\mapsto (\lambda \mu, x,a,l)
$$
identifies the latter group with the quotient
$\left(\Cal{G} (\bold{R}^{2N},
\psi )\times \Cal{G}( F\times \widehat{F},\psi_0)\right)/ \{(\lambda ,\lambda ^{-1})\}$
(where the subgroup is embedded in the center in an obvious way).
The Heisenberg representation $\Cal{H}$ of $\Cal{G} (\bold{R}^{2N}\times F\times \widehat{F},
\psi \psi_0)$ is thus identified with
the tensor product of the Heisenberg representations $\Cal{H}_1\otimes
\Cal{H}_2$ of the
two factors. We will take the Model I for $\Cal{H}_1$.
For $\Cal{H}_2$, take the space of complex functions on $F$ 
with the scalar product $\langle \phi ,\chi \rangle := \sum_{a\in F}
\phi (a) \overline{\chi (a)}$ and the action of 
$\Cal{G}( F\times \widehat{F},\psi_0)$ 
$$
(U_{(\lambda; a,l)} \phi ) (b) := \lambda l(b)\,\phi (a+b) .
\eqno(3.26)
$$
Let $\delta_a \in \Cal{H}_2$ be the delta function
supported by $a\in F,$  and consider the vectors $f_{T,a}:=f_T\otimes \delta_a$
in $\Cal{H}$. 

\smallskip

Fix a lattice embedding  $D\subset \bold{R}^{2N}\times F\times \widehat{F}$
and denote by $D_0$ the kernel of the projection $D\to F\times \widehat{F}$.
Since $\psi_0$ is not antisymmetric, we will assume that a map
$D\to\bold{C}_1^*,\,h\mapsto c_h$ has been chosen
such that $c_h$ depends only on $h+D_0$, is 1 on $D_0$, 
and the form on $D\times D$
$$
\alpha (g,h) := \frac{c_g c_h}{c_{g+h}}\,
\psi (g,h) \psi_0 (g,h)
\eqno(3.27)
$$
is antisymmetric. 
\smallskip

Consider $\Cal{H}$ as a $C^*(D,\alpha )$--module via
$e_{D,\alpha} (h) \mapsto E(h) := c_h \,U_{(1,h)}$.

\medskip

{\bf 3.6.1. Theorem.} {\it (i) The scalar products 
${}_D\langle f_{T,a} , f_{T,b}\rangle$
are quantum theta functions belonging to
the space $\Gamma (\Cal{L} )$ where $\Cal{L}$
is the multiplier
$$
D_0 \to \Cal{G} (D,\alpha ): \  g\mapsto [C_g; x_g, g],
\eqno(3.28)
$$
$C_g, x_g$ being defined by (3.10), with $H$ is lifted to
$\bold{R}^N\times F\times \widehat{F}$ via projection.

\smallskip

(ii) The scalar products ${}_D\langle f_{T,a} , f_{T,b}\rangle$ 
form a basis of $\Gamma (\Cal{L} )$.
}

\medskip

{\bf Proof.} We have for $h = (h^{\prime}, a_h, l_h)\in D,\,h^{\prime}\in 
\bold{R}^N,\, a_h\in F,\, l_h\in \widehat{F},$
$$
e_{D,\alpha}(h)\,f_{T,b} = c_h U_{(1,h)}\,(f_T\otimes \delta_b) =
$$
$$
c_h U_{(1,h^{\prime})}f_T \otimes U_{(1,a_h,l_h)}\delta_b =
c_h U_{(1,h^{\prime})}f_T \otimes l_h\delta_{b-a_h} .
\eqno(3.29)
$$ 
Therefore, 
$$
\langle f_{T,a}, e_{D,\alpha}(h)\,f_{T,b} \rangle =
\langle f_T, c_h U_{(1,h^{\prime})}f_T \rangle \cdot
\langle \delta_a , l_h\delta_{b-a_h}\rangle =
$$
$$
\frac{\overline{c}_h\overline{l}_h(a)\,\delta_{a+a_h,b}}{\sqrt{2^N\,\roman{det\,Im}\,T}}\,
 e^{-\frac{\pi}{2}\,H(\underline{h}^{\prime},\underline{h}^{\prime})} 
\eqno(3.30)
$$
(cf. (3.8) and (3.26)). Hence
$$
{}_D\langle f_{T,a} , f_{T,b}\rangle =
\frac{1}{\sqrt{2^N\,\roman{det\,Im}\,T}}
\sum_{h\in D} \overline{c}_h\overline{l}_h(a)\,\delta_{a+a_h,b}\,e^{-\frac{\pi}{2}\,
H(\underline{h}^{\prime},\underline{h}^{\prime})}\,
e_{D,\alpha}(h).
\eqno(3.31)
$$
Denote the last sum  $\Theta_{a,b}$, and take  $g\in D_0$. We have
$$
C_ge_{D,\alpha}(g)\,x_g^*(\Theta_{a,b})=
$$
$$
e^{-\frac{\pi}{2}H(\underline{g}^{\prime},\underline{g}^{\prime})}\, e_{D,\alpha}(g)
\, \sum_{h\in D} \overline{c}_h\overline{l}_h(a)\,\delta_{a+a_h,b}\,
e^{-\frac{\pi}{2}\,H(\underline{h}^{\prime},\underline{h}^{\prime})
-\pi H(\underline{g}^{\prime},\underline{h}^{\prime})}\,
e_{D,\alpha}(h)=
$$
$$
\sum_{h\in D} \overline{c}_{g+h}\overline{l}_{g+h}(a)\,\delta_{a+a_{g+h},b}\,
e^{-\frac{\pi}{2}\,H(\underline{g}^{\prime}+\underline{h}^{\prime},
\underline{g}^{\prime}+\underline{h}^{\prime})}\,
e_{D,\alpha}(g+h) =\Theta_{a,b}
$$
because in view of (3.27), $e_{D,\alpha}(g)\,e_{D,\alpha}(h)=
e^{\pi i\,\roman{Im}\,H(\underline{g}^{\prime},\underline{h}^{\prime})}\,e_{D,\alpha}(g+h)$
whenever $g\in D_0,\,h\in D,$
and moreover, $\overline{c}_h\overline{l}_h(a)\,\delta_{a+a_h,b}$ depends only
on $h\,\roman{mod}\,D_0.$

\medskip

(ii) From [Ma3], Theorem 2.4.1, it follows
that $\roman{dim}\,\Gamma (\Cal{L})= [D:D_0]$. The latter index
equals $\roman{card}\,F\times\widehat{F} = (\roman{card}\,F)^2$
because $D\subset \bold{R}^N\times  F\times\widehat{F}$
is a lattice embedding and hence the map
$h\mapsto (a_h,l_h)$ must be surjective. On the other hand,
when $a,b$ run over $F$, the functions
$\phi_{a,b}:\,D/D_0 \to \bold{C}:\, h +D_0\mapsto \overline{c}_hl_h(a)\,\delta_{a+a_h,b}$
span he whole space of functions. From this one
derives the last statement of the Theorem.

\medskip

{\bf 3.7. Further problems.} The picture described above is incomplete in at least
two respects.

\smallskip

First, we did not treat theta vectors in all possible lattice embeddings,
namely embeddings into self--dual locally compact groups of the form
{\it vector space} $\times$ {\it finite group} $\times$ {\it torus} $\times$ {\it lattice}.
Besides the naturality of this question, it is necessary to understand
the situation because the Heisenberg modules produced
from such embeddings provide some useful generators of
the strong Morita equivalence group $SO(n,n;\bold{Z})$, see [Li].
More generally, one can try to treat directly the scalar products
of theta vectors in the toric projective modules endowed with
an Hermitean connection of constant curvature. Are they all
quantum thetas? Do we get in this way quantum
thetas more general than those described in the Theorems 3.5.1 
and 3.6.1? 

\smallskip

Second, one should systematically study the behavior of theta vectors 
and their scalar products
with respect to the tensor products of toric bimodules,
as this was done for two--dimensional tori in [PoS].   

\newpage

\centerline{\bf References}

\medskip

[Bo] F.~Boca. {\it Projections in rotation algebras and theta functions.}
Comm. Math. Phys., 202 (1999), 325--357.

\smallskip

[DiSch] M.~Dieng, A.~Schwarz. {\it Differential
and
complex geometry of two--dimensional noncommutative tori.}
e--Print  math.QA/0203160

\smallskip

[Li] Hanfeng Li. {\it Strong Morita equivalence of higher
dimensional noncommutative tori.} e--Print math.OA/0303123

\smallskip

[Ma1] Yu.~Manin. {\it  Quantized theta--functions.} In: Common
Trends in Mathematics and Quantum Field Theories (Kyoto, 1990), 
Progress of Theor. Phys. Supplement, 102 (1990), 219--228.

\smallskip

[Ma2] Yu.~Manin. {\it Mirror symmetry and quantization of abelian varieties.}
In: Moduli of Abelian Varieties, ed. by C.~Faber et al.,
Progress in Math., vol. 195, Birkh\"auser, 2001, 231--254.
e--Print  math.AG/0005143

\smallskip

[Ma3] Yu.~Manin. {\it Theta functions, quantum tori and Heisenberg groups}.
Lett. in Math. Physics, 56:3 (2001)
(special issue, Euroconference M.~Flato, Part III), 295--320.
e--Print math.AG/0011197

\smallskip

[Ma4] Yu.~Manin. {\it  Real multiplication and noncommutative
geometry.} e--Print math.AG/0202109

\smallskip

[Mu] D.~Mumford (with M.~Nori and P.~Norman). 
{\it Tata Lectures on Theta III.} Progress in Math., vol.~97,
Birkh\"auser, 1991.

\smallskip

[PoS] A.~Polishchuk, A.~Schwarz. {\it Categories of holomorphic vector
bundles on noncommutative two--tori.} e--Print math.QA/0211262

\smallskip

[Ri1] M.~A.~Rieffel. {\it Strong Morita equivalence of
certain transformation group $C^*$--algebras.}
Math. Annalen, 222 (1976), 7--23.

\smallskip

[Ri2] M.~A.~Rieffel. {\it Von Neumann algebras associated with pairs of
lattices in Lie groups.} Math.~Ann., 257 (1981), 403--418.

\smallskip

[Ri3] M.~A.~Rieffel. {\it $C^*$--algebras associated with irrational rotations.}
Pacific J.~Math., 93 (1981), 415--429.

\smallskip

[Ri4] M.~A.~Rieffel. {\it The cancellation theorem for projective
modules over irrational rotation $C^*$--algebras.}
Proc.~Lond.~Math.~Soc. (3), 47 (1983), 285--303.

\smallskip

[Ri5] M.~A.~Rieffel. {\it Projective modules over higher--dimensional
non--commutative tori.} Can.~J.~Math., vol.~XL, No.~2 (1988), 257--338.

\smallskip

[Ri6] M.~A.~Rieffel. {\it Non--commutative tori --- a case
study of non--commutative differential manifolds.}
In: Cont.~Math., 105 (1990), 191--211.

\smallskip

[RiSch] M.~A.~Rieffel, A.~Schwarz. {\it Morita equivalence
of multidimensional non--commutative tori.} 
Int. J. Math., 10 (1999), 289--299. e--Print   math.QA/9803057

\smallskip

[Sch] A.~Schwarz. {\it Theta--functions on non--commutative tori.}
e--Print \newline math.QA/0107186

\bigskip

e-mail: manin\@mpim-bonn.mpg.de

\enddocument